\documentclass[11]{article}
\usepackage{amssymb}
\input{epsf.sty}
\newtheorem{theorem}{Theorem}
\newtheorem{proposition}{Proposition} 

\newcommand{\sss}{\scriptscriptstyle}

\newcommand{\sgn}{\mathop{\rm sgn}}

\newcommand{\ds}{\displaystyle}

\newcommand{\barF}{\overline{F}}

\newcommand{\Prob}{\hbox{${\bf P}$}}

\newcommand{\fine}{\hfill $\Box$}         % box vuoto
\begin{document}
\date{\empty}
\newenvironment{proof}{\begin{trivlist}
\item[\hspace{\labelsep}{\bf\noindent Proof. }]}
{\end{trivlist}}
\title{
\Large\bf
Random motion with gamma-distributed\\
alternating velocities in biological modeling\thanks{Paper appeared in: \newline
R.\  Moreno-D\'{\i}az et al.\ (Eds.):    EUROCAST 2007,  
Lecture Notes in Computer Science, Vol.\ 4739, pp.\ 163--170, 2007. 
Springer-Verlag, Berlin, Heidelberg. ISBN: 978-3-540-75866-2.}
}
\author{\bf Antonio Di Crescenzo \rm and \bf Barbara Martinucci\\
\hfill\\
\rm Dipartimento di Matematica e Informatica, Universit\`a di Salerno\\
Via Ponte don Melillo,  I-84084 Fisciano (SA), Italy\\
Email: \{adicrescenzo,bmartinucci\}@unisa.it 
}
\pagestyle{plain}
\maketitle
\begin{abstract}
Motivated by applications in mathematical biology concerning randomly  
alternating motion of micro-organisms, we analyze a generalized 
integrated telegraph process. The random times between consecutive 
velocity reversals are gamma-distributed, and perform an alternating 
renewal process. We obtain the probability law and the mean 
of the process. 
\end{abstract}
%
%--------------------------------------------------------------------
\section{Introduction}
\label{sec1}
%--------------------------------------------------------------------
%
The telegraph random process describes the motion of a particle on the real line, 
traveling at constant speed, whose direction is reversed at the arrival epochs 
of a Poisson process. After some initial works, such as \cite{Go51}, \cite{Ka74} 
and \cite{Or90a}, numerous efforts have been made by numerous authors and through 
different methods to analyze this process. Various results on the telegraph 
process, including the first-passage-time density and the distribution of motion 
in the presence of reflecting and absorbing barriers have been obtained 
in \cite{Fo92}, \cite{FoKa94} and \cite{Or95}. A wide and comprehensive review 
devoted to this process has recently been offered by Weiss \cite{Weiss02}, 
who also emphasized its relations with some physical problems. 
\par
In various applications in biomathematics the telegraph process arises as a 
stochastic model for systems driven by dichotomous noise (see \cite{DiCrMa2006}, 
for instance). Two stochastic processes modeling the major 
modes of dispersal of cells or organisms in nature are introduced in \cite{Othmer}; 
under certain assumptions, the motion consisting of sequences of runs 
separated by reorientations with new velocities is shown to be governed 
by the telegraph equation. Moreover, the discrete analog of 
the telegraph process, i.e.\ the correlated random walk, is usually 
used as a model of the swarming behavior of myxobacteria 
(see \cite{Erdetal2004}, \cite{Hill}, \cite{Komin} and \cite{Lutscher}). 
Processes governed by hyperbolic equations are also used to describe movement 
and interaction of animals \cite{Lu2002} and chemotaxis \cite{HiSt2000}.
Moreover, the integrated telegraph process has been also used to model 
wear processes \cite{DiPe2000} and to describe the dynamics of the price 
of risky assets \cite{DiPe2002}. 
\par
Many authors proposed suitable generalizations of the telegraph process, 
such as the 1-dimensional cases with three cyclical velocities \cite{Or90b}, 
or with $n$ values of the velocity \cite{Ko98}, or with random 
velocities \cite{StZa04}. See also the paper by Lachal \cite{Latchal}, where 
the cyclic random motion in $\mathbb{R}^d$ with $n$ directions is studied.  
\par
A generalized integrated telegraph process whose random times separating 
consecutive velocity reversals have a general distribution and perform an 
alternating renewal process has been studied in \cite{DiCrescenzo2001} and \cite{Za03}. 
Along the line of such articles, in this paper we study a stochastic model 
for particles motion on the real line with two alternating velocities $c$ 
and $-v$. The random times between consecutive reversals of direction perform 
an alternating renewal process and are gamma distributed, which extends 
the Erlang-distribution case treated in \cite{DiCrescenzo2001}. 
\par
In Section~\ref{sec2} we introduce the stochastic process $\{(X_t,V_t); t\geq 0\}$, 
with $X_t$ and $V_t$ denoting respectively position and velocity of the particle 
at time $t$. In Section~\ref{sec3} we obtain a series-form of the random motion 
probability law for gamma-distributed random inter-renewal times, whereas the 
mean value of $X_t$ conditional on initial velocity is finally obtained in 
Section~\ref{sec4}.  
%
%----------------------------------------------------------------------
\section{The random motion}
\label{sec2}
%----------------------------------------------------------------------
We consider a random motion on $\mathbb{R}$ with two alternating velocities 
$c$ and $-v$, with $c,v>0$. The direction of motion is forward or backward 
when the velocity is $c$ or $-v$, respectively. Velocities change according 
to the alternating counting process $\{N_t; t\geq 0\}$ characterized by 
renewal times $T_1,T_2,\ldots$, so that $T_n$ is the $n$-th random instant 
in which the motion changes velocity. Hence, 
$$
 N_0=0, \qquad N_t=\sum_{n=1}^{\infty}{\bf 1}_{\{T_n\leqslant t\}},\quad t>0.
$$ 
Let $\{(X_t,V_t); t\geq 0\}$ be a stochastic process on $\mathbb{R}\times\{-v,c\}$, 
where $X_t$ and $V_t$ give respectively position and velocity of the motion at 
time $t$. Assuming that $X_0=0$ and $v_0\in\{-v,c\}$, for $t>0$ we have:  
\begin{equation} 
 X_t =\int_0^t V_s \,{\rm d}s, 
 \qquad
 V_t=\frac{1}{2}(c-v)+\sgn(V_0)\,\frac{1}{2}(c+v)\,(-1)^{N_t}.
 \label{eq:36}
\end{equation}
Denoting by $U_k$ ($D_k$) the duration of the $k$-th time interval during 
which the motion goes forward (backward), we assume that $\{U_k; k=1,2,\ldots\}$ 
and $\{D_k; k=1,2,\ldots\}$ are mutually independent sequences of independent copies 
of non-negative and absolutely continuous random variables $U$ and $D$. 
\par
If the motion does not change velocity in $[0,t]$, then $X_t=V_0\,t$. 
Otherwise, if there is at least one velocity change in $[0,t]$, then 
$-vt<X_t<ct$ w.p.\ 1. Hence, the conditional law of 
$\{(X_t,V_t);t\geq 0\}$ is characterized by a discrete component 
$$
\Prob\{X_t=yt, V_t=y \,|\,X_0=0, V_0=y\}, 
$$
and by an absolutely continuous component 
\begin{equation}
 p(x,t\,|\,y)= f(x,t\,|\,y) + b(x,t\,|\,y),
 \label{equation:p}
\end{equation} 
where
$$
 f(x,t\,|\,y) = {\partial \over \partial x}\Prob\{X_t\leq x,V_t=c\,|\,X_0=0, V_0=y\}, 
$$
$$
 b(x,t\,|\,y) = {\partial \over \partial x}\Prob\{X_t\leq x,V_t=-v\,|\,X_0=0, V_0=y\}, 
$$
with $t>0$, $-vt<x<ct$ and $y\in\{-v,c\}$. 
\par
The formal conditional law of $\{(X_t,V_t);t\geq 0\}$ has been given in Theorem 2.1 
of \cite{DiCrescenzo2001} for $V_0=c$. Case $V_0=-v$ can be treated by 
symmetry. 
\par
Explicit results for the probability law have been obtained in Theorem 3.1 of 
\cite{DiCrescenzo2001} when the random times $U$ and $D$ separating consecutive 
velocity reversals have Erlang distribution. This case describes the random 
motion of particles subject to collisions arriving according to a Poisson 
process with rate $\lambda$ if the motion is forward and rate $\mu$ it is 
backward. When the motion has initial velocity $c$ ($-v$), then the first $n-1$ 
($r-1$) collisions have no effect, whereas the $n$th ($r$th) collision causes 
a velocity reversal. In the following section we shall treat the more general 
case in which the random inter-renewal times are gamma distributed.
%---------------------------------------------------------------
\section{\bf Gamma-distributed random times}
\label{sec3}
%---------------------------------------------------------------
We assume that the random times $U$ and $D$ are gamma distributed with 
parameters ($\lambda$,$\alpha$) and ($\mu$,$\beta$), respectively, where 
$\lambda, \mu>0$ and $\alpha, \beta>0$. Hereafter we obtain the probability 
law of $\{(X_t,V_t); t\geq 0\}$ for this case.
\begin{theorem}
If $U$ and $D$ are gamma-distributed with parameters $(\lambda,\alpha)$ 
and $(\mu,\beta)$, respectively, for $t>0$ it is
\begin{equation}
 \Prob\{X_t=ct, V_t=c\,|\,X_0=0, V_0=c\}=
\frac{\Gamma(\alpha,\lambda t)}{\Gamma(\alpha)},
\label{discreta}
\end{equation}
and, for $-vt<x<ct$, 
\begin{equation}
\hspace*{-0.2cm} 
f(x,t\,|\,c)=\frac{1}{c+v}\left\{e^{-\mu\overline{x}}
\sum_{k=1}^{+\infty}\frac{\mu^{k\beta}(\overline{x})^{k\beta-1}}{\Gamma(k\beta)}\bigg[P(k\alpha,\lambda x^*)
-P(k\alpha+\alpha,\lambda x^*)\bigg]\right\},
\label{gammaf}
\end{equation}
\begin{eqnarray}
&& \hspace*{-1.cm} 
b(x,t\,|\,c)=\frac{1}{c+v}\Biggl\{\frac{\lambda^{\alpha}e^{-\lambda x^*}
(x^*)^{\alpha-1}}{\Gamma(\alpha)}\frac{\Gamma(\beta,\mu\overline{x})}{\Gamma(\beta)} 
\nonumber \\
&& \hspace*{0.4cm} 
+ e^{-\lambda x^*} \sum_{k=1}^{+\infty}\frac{\lambda^{(k+1)\alpha}(x^*)^{(k+1)\alpha-1}}
{\Gamma((k+1)\alpha)} 
\biggl[P(k\beta,\mu\overline{x})-P(k\beta+\beta,\mu\overline{x})\biggr]\Biggr\},
\label{gammab} 
\end{eqnarray}
where 
$$
 \overline{x}=\overline{x}(x,t)=\frac{ct-x}{c+v}, 
 \qquad 
 x^{*}=x^{*}(x,t)=\displaystyle\frac{vt+x}{c+v},
$$
and 
\begin{equation}
 \Gamma(a,u)=\int_u^{\infty} t^{a-1}e^{-t} {\rm d}t,
 \qquad 
 P(a,u)=\frac{1}{\Gamma(a)}\int_0^{u} t^{a-1}e^{-t} {\rm d}t,
 \qquad a>0.
\label{gammaP}
\end{equation}
\label{teorema3}
\end{theorem}
\begin{proof}
Making use of (2.4) of \cite{DiCrescenzo2001} and noting that for $k\geq 1$ 
the pdfs of $U^{(k)}=U_1+\ldots+U_k$ e $D^{(k)}=D_1+\ldots+D_k$ are given by
\begin{eqnarray}
f_U^{(k)}(x)=\frac{\lambda^{k\alpha}x^{k\alpha-1}e^{-\lambda x}}{\Gamma(k\alpha)},
\qquad f_D^{(k)}(x)=\frac{\mu^{k\beta}x^{k\beta-1}e^{-\mu x}}{\Gamma(k\beta)}, 
\qquad x>0, 
\label{eq:319} 
\end{eqnarray}
we have
\begin{equation}
f(x,t\,|\,c)=\frac{1}{c+v}\,e^{-\mu\overline{x}}
e^{\lambda\overline{x}}\sum_{k=1}^{+\infty}\frac{\mu^{k\beta}(\overline{x})^{k\beta-1}}{\Gamma(k\beta)} 
\frac{\lambda^{k\alpha}}{\Gamma(k\alpha)\Gamma(\alpha)} \, {\cal I}_k,
\label{gammaf2}
\end{equation}
where 
\begin{equation}
{\cal I}_k:=\int_{\overline{x}}^t e^{-\lambda s}
(s-\overline{x})^{k\alpha-1}\Gamma(\alpha,\lambda(t-s))\,{\rm d}s, 
 \qquad k\geq 1.
\label{integI}
\end{equation}
Noting that, due to (\ref{gammaP}), $\Gamma(a,u)=\Gamma(a)\,\big[1-P(a,u)\big]$ 
we obtain
\begin{equation}
{\cal I}_k={\cal I}_{1,k}-{\cal I}_{2,k},
\label{integI1I2}
\end{equation}
where, for $k\geq 1$  
\begin{equation}
 {\cal I}_{1,k}:=
 \Gamma(\alpha)\int_{\overline{x}}^t e^{-\lambda s}
 (s-\overline{x})^{k\alpha-1}\,{\rm d}s
 =\Gamma(k\alpha)\,\Gamma(\alpha)\, 
 e^{-\lambda\overline{x}}\lambda^{-k\alpha} P(k\alpha,\lambda x^*), 
\label{eq:323}
\end{equation}
\begin{eqnarray}
 && \hspace{-1.2cm} 
 {\cal I}_{2,k}:= \Gamma(\alpha)\int_{\overline{x}}^t e^{-\lambda s}
 (s-\overline{x})^{k\alpha-1}\,P(\alpha,\lambda(t-s))\,{\rm d}s 
 \nonumber \\
 && \hspace{-0.4cm} 
  = e^{-\lambda\overline{x}}\Gamma(\alpha)\int_0^{x^*} e^{-\lambda y}
 y^{k\alpha-1}\,P(\alpha,\lambda(x^*-y))\,{\rm d}y
 = e^{-\lambda\overline{x}} \lambda^{-k\alpha} \,G(\lambda x^*),
 \label{eq:324}
\end{eqnarray}
with
\begin{equation}
 G(\lambda x^*):=\int_0^{\lambda x^*} e^{-u}u^{k\alpha-1}
 \Biggl(\int_0^{\lambda x^*-u}e^{-\tau}\tau^{\alpha-1} \,{\rm d}\tau\Biggr){\rm d}u.
 \label{defG}
\end{equation}
Making use of the Laplace transform  of $(\ref{defG})$ it follows 
$$
 {\cal L}_z\{G(\lambda x^*)\}
 ={\cal L}_z\left\{ e^{-{\lambda x^*}}\left(\lambda x^*\right)^{k\alpha-1}\right\}
{\cal L}_z\Biggl\{\int_0^{\lambda x^*} e^{-\tau} \tau^{\alpha-1}{\rm d}\tau\Biggr\}
=\frac{\Gamma(k\alpha)\,\Gamma(\alpha)}{z(z+1)^{k\alpha+\alpha}}.
$$
Hence, from identity 
$$
 {\cal L}_z\left\{P(k\alpha+\alpha,\lambda x^*)\right\}
 ={\cal L}_z\left\{\int_0^{\lambda x^*} {u^{k\alpha+\alpha-1}e^{-u}\over \Gamma(k\alpha+\alpha)}\,{\rm d}u\right\}
 =\frac{1}{z(z+1)^{k\alpha+\alpha}},
$$
we have 
$$
 G(\lambda x^*)= \Gamma(k\alpha)\,\Gamma(\alpha)\,P(k\alpha+\alpha,\lambda x^*).
$$
Eqs.\ (\ref{integI1I2})$\div$(\ref{defG}) thus give  
\begin{equation}
 {\cal I}_k
 =\Gamma(k\alpha)\,\Gamma(\alpha)\, e^{-\lambda\overline{x}}
 \lambda^{-k\alpha}\,\big[P(k\alpha,\lambda x^*) - P(k\alpha+\alpha,\lambda x^*)\big].
 \label{valIk}
\end{equation}
Eq. (\ref{gammaf}) then follows from (\ref{gammaf2}) and (\ref{defG}). 
In order to obtain $b(x,t\,|\,c)$, we recall (2.5) of \cite{DiCrescenzo2001} 
and make use of (\ref{eq:319}) to obtain 
\begin{eqnarray}
b(x,t\,|\,c)&=& \frac{1}{c+v}\bigg\{\frac{\lambda^{\alpha}e^{-\lambda x^*}
(x^*)^{\alpha-1}}{\Gamma(\alpha)}\frac{\Gamma(\beta,\mu\overline{x})}{\Gamma(\beta)} 
\nonumber \\
&+& \left. e^{-\lambda x^*}e^{\mu x^*}\sum_{k=1}^{+\infty}\frac{\lambda^{(k+1)\alpha}
(x^*)^{(k+1)\alpha-1}}{\Gamma((k+1)\alpha)} \right.
\nonumber \\
&\times&   \frac{\mu^{k\beta}}{\Gamma(\beta)\Gamma(k\beta)}
\int_{x^*}^t e^{-\mu s}(s-x^*)^{k\beta-1} \Gamma(\beta,\mu(t-s)){\rm d}s\bigg\}. 
\label{eq:328}
\end{eqnarray}
Due to (\ref{integI}), the integral in (\ref{eq:328}) can be calculated from 
(\ref{valIk}) by interchanging $x^*$, $\beta$, $\mu$ with $\overline{x}$, $\alpha$,  
$\lambda$, respectively. Eq.\ (\ref{gammab}) then follows after some calculations.
\fine
\end{proof}
%
%%%%%%%%%%%%%%%%%%%%%%%%%%%%%%%%%%%%%%%%%%%%%%%%%%%%%%%%%%%%%%%%%%%%%%%%
\begin{figure}[t]
\begin{center}
%\hspace*{-0.3cm}
\epsfxsize=9cm
\centerline{\epsfbox{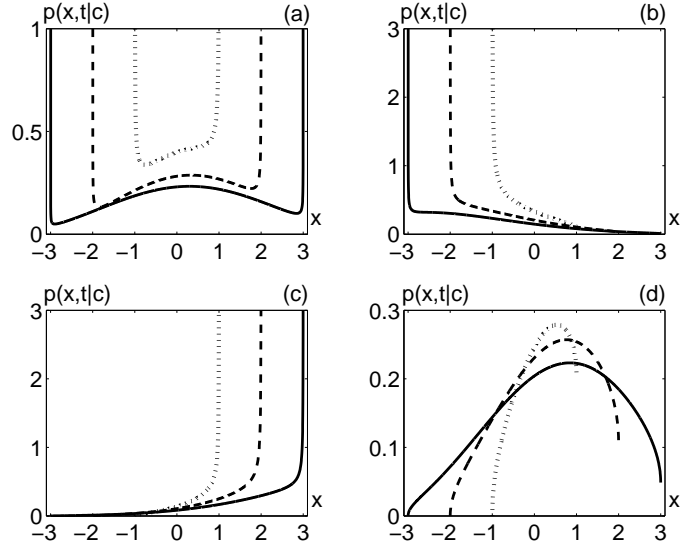}}
\end{center}
\vspace{-0.5cm}
\caption{Density (\ref{equation:p}) for $c=v=1$, $\lambda=\mu=1$, and $t=1$ (dotted line), 
$t=2$ (dashed line), and $t=3$ (solid line), with (a) $\alpha=\beta=0.5$, 
(b) $\alpha=0.5$ and $\beta=1.5$, (c) $\alpha=1.5$ and $\beta=0.5$, (d) $\alpha=\beta=1.5$.
}
\end{figure}
%%%%%%%%%%%%%%%%%%%%%%%%%%%%%%%%%%%%%%%%%%%%%%%%%%%%%%%%%%%%%%%%%%%%%%%%
%
\par 
Figure 1 shows density $p(x,t\,|\,c)$ as $x$ varies in $(-vt,ct)$ for various 
choices of $t$, $\alpha$ and $\beta$. 
Hereafter we analyze the obtain the limits of densities $f(x,t\,|\,c)$ and 
$b(x,t\,|\,c)$ at the extreme points of interval $(-vt, ct)$, for any fixed $t$. 
\begin{proposition}\label{proposizione1}
Under the assumptions of Theorem \ref{teorema3} we have 
$$
\lim_{x \downarrow -vt} f(x,t|c)=0, 
\quad\!
\lim_{x \uparrow ct} f(x,t|c)
 = \cases{+\infty, &  $0<\beta<1$ \cr
\ds\frac{\mu}{c+v}\big[P(\alpha,\lambda t)-P(2 \alpha, \lambda t)\big], 
&  $\beta=1$ \cr
0, &  $\beta>1$,}
$$
$$
\lim_{x \uparrow ct} b(x,t|c)
 =\frac{\lambda^{\alpha} {\rm e}^{-\lambda t} t^{\alpha-1}}{(c+v) \Gamma(\alpha)},
\quad\!
\lim_{x \downarrow -vt} b(x,t|c)
 = \cases{+\infty, & $0<\alpha<1$\cr
\ds\frac{\lambda\,\Gamma(\beta, \mu t)}{(c+v) \Gamma(\alpha) \Gamma (\beta)}, 
& $\alpha=1$ \cr
0, & $\alpha>1$. \cr}
$$
\end{proposition}
\par
From Proposition \ref{proposizione1} we note that if $\alpha<1$ $(\beta<1)$, 
i.e.\ the gamma inter-renewal density has a decreasing hazard rate, then the 
backward (forward) density is divergent when $x$ approaches $-vt$ $(ct)$. 
This is very different from the behavior exhibited in the case of 
Erlang-distributed inter-renewals (see Corollary $3.1$ 
of \cite{DiCrescenzo2001}), when the limits are finite. 
%
%-------------------------------------------------
\section{\bf Mean value}
\label{sec4}
%-------------------------------------------------
In this Section we obtain the mean value of $X_t$   
when random times $U$ and $D$ are identically gamma distributed. 
\begin{theorem}\label{teorema4}
Let $U$ and $D$ have gamma distribution with parameters $(\lambda, \alpha)$. 
For any fixed $t\in(0,+\infty)$, we have 
\begin{equation}
 E\big[X_t\,\big|\,V_0\big]
 = V_0\,t+ \frac{c+v}{\lambda} \sgn(V_0)
 \sum_{k=1}^{+\infty} (-1)^k \Big[\lambda t \, P(k \alpha, \lambda t)
 - k \alpha P(k \alpha+1, \lambda t) \Big].
\label{media_}
\end{equation}
\end{theorem}
\begin{proof}
Due to Eqs.\ (\ref{eq:36}) and recalling that $\Prob(T_{k}\leq s)=P(k\alpha,\lambda s)$, 
$s\geq 0$, it is
\begin{eqnarray}
 E\big[X_t\,\big|\,V_0\big] \!\!\!
 &=& \!\!\! \frac{1}{2}(c-v)t
 + \frac{1}{2}(c+v)\sgn(V_0)\int_0^t E\left[(-1)^{N_s}\right]{\rm d}s  
 \label{medsum} \\
 &=& \!\!\! \frac{1}{2}(c-v)t
 + \frac{1}{2}(c+v)\sgn(V_0)\int_0^t 
 \left\{1+2\sum_{k=1}^{+\infty} (-1)^k \,P(k\alpha,\lambda s)\right\}{\rm d}s 
 \nonumber \\
 &=& \!\!\! V_0\,t +(c+v)\sgn(V_0)\sum_{k=1}^{+\infty} (-1)^k 
 \int_0^t P(k\alpha,\lambda s)\,{\rm d}s.  
 \nonumber
\end{eqnarray}
(Note that the above series is uniformly convergent.) 
Moreover, recalling (\ref{gammaP}) it is not hard to see that 
$$
 \int_0^t P(k\alpha,\lambda s)\,{\rm d}s 
 = t\,P(k\alpha,\lambda t)-{k \alpha \over \lambda}\,P(k \alpha +1,\lambda t).
$$
Eq.\ (\ref{media_}) then immediately follows. 
\fine
\end{proof}
%
%%%%%%%%%%%%%%%%%%%%%%%%%%%%%%%%%%%%%%%%%%%%%%%%%%%%%%%%%%%%%%%%%%%%%%%%
%
\begin{figure}[t]
\begin{center}
%\hspace*{-0.3cm}
\epsfxsize=8cm
\centerline{\epsfbox{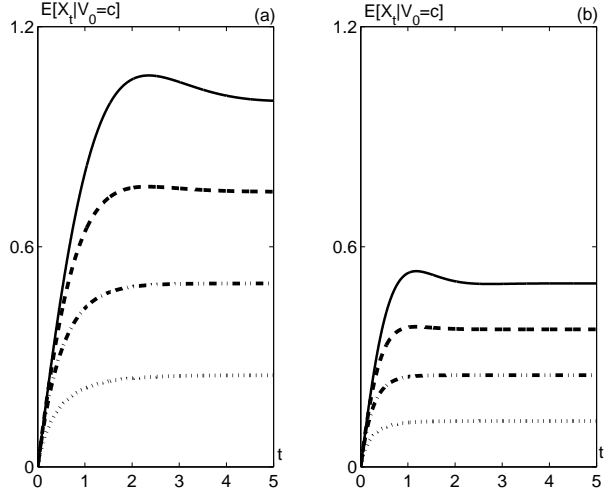}}
\end{center}
\vspace{-0.5cm}
\caption{Mean value $E[X_t\,|\,V_0=c]$, for $c=v=1$ and $\alpha=0.5$ (dotted line), 
$\alpha=1$ (dash-dot line), $\alpha=1.5$ (dash line), $\alpha=2$ (solid line), 
with (a) $\lambda=1$ and (b) $\lambda=2$.} 
\label{figmedia}
\end{figure}
%%%%%%%%%%%%%%%%%%%%%%%%%%%%%%%%%%%%%%%%%%%%%%%%%%%%%%%%%%%%%%%%%%%%%%%%
\par
The graphs given in Figure \ref{figmedia} show the mean value of $X_t$ 
conditional on $V_0=c$ for some choice of the involved parameters. 
We note that, being 
$P(\alpha,t) \sim {t^{\alpha-1}/\Gamma(\alpha)}$ as $t\to 0$, 
under the assumptions of Theorem \ref{teorema4} from (\ref{media_}) we have 
$$
 E\big[X_t\,\big|\,V_0\big]\sim V_0\,t  
 \qquad 
 \hbox{as }t\to 0. 
$$
\par
We remark that when $\alpha=n$ is integer, i.e.\ the random times $U$ and $D$ 
are Erlang-distributed with parameters $(\lambda, n)$, then 
$E\big[X_t\,\big|\,V_0\big]$ can be computed making use of 
(\ref{medsum}) and noting that  
$$
 E\big[(-1)^{N_s}\big]
 = 1-2{\rm e}^{-\lambda s}\sum_{k=0}^{+\infty}\, 
 \sum_{j=2nk+n}^{2nk+2n-1}{(\lambda s)^j\over j!}.  
$$
For instance, in this case the following expressions hold for $t>0$: 
\begin{center}
\begin{tabular}{ll}
\hline
$n$ & $E\big[X_t\,\big|\,V_0\big]$ \\
\hline
$1\quad$ & 
$\ds\frac{(c-v)t}{2}+\frac{(c+v)}{4 \lambda}\sgn(V_0)\,[1-{\rm e}^{-2 \lambda t}]$ \\
$2\quad$ &
$\ds\frac{(c-v)t}{2}+\frac{(c+v)}{2 \lambda}\sgn(V_0)\,[1-{\rm e}^{- \lambda t} \cos(\lambda t)]$ \\
$3\quad$ &
$\ds\frac{(c-v)t}{2}+\frac{(c+v)}{2 \lambda}\sgn(V_0)\left\{\frac{1-{\rm e}^{-2 \lambda t}}{6}
+\frac{4}{3}[1-{\rm e}^{-\frac{\lambda t}{2}}\cos(\frac{\sqrt{3}}{2}\lambda t)]\right\}$ \\
$4\quad$ &
$\ds\frac{(c-v)t}{2}+\frac{(c+v)}{2 \lambda}\sgn(V_0)\,\left\{
\big[1-(1+\frac{\sqrt{2}}{2}) {\rm e}^{-\lambda t (1-\frac{\sqrt{2}}{2})}
\cos(\frac{\sqrt{2}}{2}\lambda t)\big] \right.$ \\
{} & 
$\left. + \big[1-(1-\ds\frac{\sqrt{2}}{2}) {\rm e}^{-\lambda t (1+\frac{\sqrt{2}}{2})}
\cos(\frac{\sqrt{2}}{2}\lambda t)\big]\right\}$ \\
\hline
\end{tabular}
\end{center}
%
%--------------------------------------------------------------------

%

\begin{thebibliography}{99}
%--------------------------------------------------------------------
%
\bibitem{DiCrescenzo2001}
Di Crescenzo, A.: On random motions with velocities alternating at Erlang-distributed 
random times. Advances in Applied Probability, {33} (2001) 690--701.
%
\bibitem{DiCrMa2006}
Di Crescenzo, A. and Martinucci, B.:
On the effect of random alternating perturbations on hazard rates. 
Scientiae Mathematicae Japonicae, 64 (2006) 381--394. 
%
\bibitem{DiPe2000}
Di Crescenzo, A. and Pellerey, F.: Stochastic comparison of wear 
processes characterized by random linear wear rates. 
2nd International Conference on Mathematical Methods in Reliability, 
Abstracts Book, Vol. 1, eds. M. Mikulin and N. Limnios, Bordeaux, (2001) 339--342.
%
\bibitem{DiPe2002}
Di Crescenzo, A. and Pellerey, F.: On Prices' evolutions based on geometric 
telegrapher's process. Applied Stochastic Models in Business and Industry, 
18:2 (2002) 171--184.
%
\bibitem{Erdetal2004}
Erdmann, U., Ebeling, W., Schimansky-Geier, L., Ordemann, A. and Moss, F.
Active Brownian particle and random walk theories of the motions of zooplankton:
Application to experiments with swarms of Daphnia. arXiv:q-bio.PE/0404018.
%
\bibitem{Fo92}
Foong, S.K.: First-passage time, maximum displacement, 
and Kac's solution of the telegrapher equation. Physical Review A, {46} (2001) 
R707--R710. 
%
\bibitem{FoKa94} 
Foong, S.K. and Kanno, S.: Properties of the telegrapher's 
random process with or without a trap. Stochastic Processes and their 
Applications, {53} (1994) 147--173. 
%
\bibitem{Go51}
Goldstein, S.: On diffusion by discontinuous movements and the 
telegraph equation. Quarterly Journal of Mechanics and Applied Mathematics, 
{4} (1951) 129--156. 
%
\bibitem{Hill}
Hill N.A. and Hader D.P.: A biased random walk model for the trajectories 
of swimming micro-organisms. Journal of Theoretical Biology, {186} (1997) 503--526.
%
\bibitem{HiSt2000}
Hillen, T. and Stevens, A.: Hyperbolic models for chemotaxis in 1-D. 
Nonlinear Analysis: Real World Applications, {1} (2000) 409--433. 
%
\bibitem{Ka74}
Kac, M.: A stochastic model related to the telegrapher's equation.  
Rocky Mountain Journal of Mathematics, {4} (1974) 497--509. 
%
\bibitem{Ko98}
Kolesnik, A.: The equations of Markovian random evolution on 
the line. Journal of Applied Probability, {35} (1998) 27--35. 
%
\bibitem{Komin}
Komin N., Erdmann U. and Schimansky-Geier L.: Random walk 
theory applied to Daphnia Motion. Fluctuation and Noise Letters, {4}, n.1 (2004) 
151--159.   
%
\bibitem{Latchal}
Lachal A.: Cyclic random motions in $\mathbb{R}^d$-space with $n$ directions. 
ESAIM Probability and Statistics, {10} (2006) 277--316. 
%
\bibitem{Lutscher}
Lutscher F. and Stevens A.: Emerging Patterns in a Hyperbolic Model 
for Locally Interacting Cell Systems. Journal of Nonlinear Science, {12} (2002) 619--640.
%
\bibitem{Lu2002}
Lutscher, F.: Modeling alignment and movement of animals and cells. 
Journal of Mathematical Biology, {45} (2002) 234--260.
%
\bibitem{Or90a}
Orsingher, E.: Probability law, flow function, maximum 
distribution of wave-governed random motions and their connections with 
Kirchoff's laws. Stochastic Processes and their Applications, {34} (1990) 49--66. 
%
\bibitem{Or90b}
Orsingher, E.: Random motions governed by third-order equations.  
Advances in Applied Probability, {22} (1990) 915--928. 
%
\bibitem{Or95}
Orsingher, E.: Motions with reflecting and absorbing barriers 
driven by the telegraph equation. Random Operators and Stochastic Equations, 
{3} (1995) 9--21. 
%
\bibitem{Othmer} 
Othmer H.G., Dunbar S. R. and Alt W.: Models of dispersal in 
biological systems. Journal of Mathematical Biology, {26}, n.3 (1988) 263--298. 
%
\bibitem{StZa04}
Stadje, W. and Zacks, S.: Telegraph processes with random 
velocities. Journal of Applied Probability, {41} (2004) 665--678.
%
\bibitem{Weiss02}
Weiss, G.H.: Some applications of persistent random walks 
and the telegrapher's equation. Physica A, {311} (2002) 381--410. 
%
\bibitem{Za03}
Zacks, S.: Generalized integrated telegrapher process and the
distribution of related stopping times. Journal of Applied Probability, {41} (2004) 
497--507.
%
\end{thebibliography}
\end{document}